%% file: ex_article.tex
\documentclass[final,hidelinks,onefignum,onetabnum]{siamart251216}

\usepackage{microtype}
\usepackage{mathtools}
\input{ex_shared}
\ifpdf
\hypersetup{
	pdftitle={Masked Symmetric Nonnegative Matrix Factorization for Community Detection in Incomplete Networks},
	pdfauthor={A. Liu, R. Gu, and R. J. Zhang}
}
\fi

\begin{document}
	
	\maketitle
	
	\begin{abstract}
		Community detection in complex networks is frequently challenged by incomplete or noisy adjacency matrices. Traditional symmetric nonnegative matrix factorization methods typically rely on zero-imputation for unobserved entries, which compromises clustering reliability. This paper proposes a Masked Symmetric Nonnegative Matrix Factorization (Masked SymNMF) framework designed to factorize partially observed networks directly. By defining a masking operator over the observed entries, the proposed model restricts the objective evaluation exclusively to valid data. To overcome the severe non-convexity inherent in the symmetric factorization, we formulate an asymmetric relaxation penalized by a regularization term. We prove the exact penalty property of this reformulated model, establishing its theoretical equivalence to the original symmetric problem under sufficient regularization. Furthermore, an alternating nonnegative least squares framework is developed, yielding tailored update rules for Multiplicative Updates, Hierarchical Alternating Least Squares, and Projected Gradient Descent algorithms. Extensive numerical experiments on synthetic datasets and real-world networks demonstrate that the proposed Masked SymNMF outperforms baseline imputation methods across varying observation densities, providing a theoretically sound and practically efficient approach for community detection in incomplete networks.
	\end{abstract}
	
	
	
	\section{Introduction}
	\label{Section1}
	\subsection{Background and Motivation}
	Complex networks, including social networks~\cite{girvan2002community,li2017utility}, biological networks~\cite{hu2018efficiently}, and sensor networks~\cite{huang2017efficient,liu2020multistep}, play a central role in modern data analysis, systems biology, social sciences, recommendation systems, and security analysis~\cite{boccaletti2006complex,newman2003structure}. Detecting community structure, which refers to subgraphs where nodes are densely connected internally but sparsely connected to other communities~\cite{girvan2002community}, helps uncover user behavior patterns in social networks, functional modules in biological systems, and propagation mechanisms in information networks.
	Community detection has emerged as an important tool in network analysis~\cite{estrada2012structure,girvan2002community}. 
	
	However, incomplete observation is pervasive in real-world networks. Due to sampling costs, observational limitations, privacy protocols, or data noise, adjacency matrices inevitably contain numerous missing entries, leading to incomplete networks~\cite{kolda2009tensor}. Examples include unvalidated interactions in PPI networks, hidden user associations in social networks, and data loss in sensor networks. 
    Under such conditions, traditional community detection techniques fail to stably extract reliable community structures, making community detection in incomplete networks a critical problem.
	
	Traditional community detection methods, such as modularity optimization~\cite{newman2004finding}, stochastic block models (SBM)~\cite{abbe2018community}, and spectral clustering~\cite{jin2021survey}, encounter bottlenecks when dealing with missing observations. While recent nonlinear relaxations have advanced modularity maximization~\cite{tudisco2018community}, and spectral methods have been adapted for sparse data~\cite{van2013community}, most existing models are highly sensitive to the edge sparsity and noise~\cite{ding2019sparsity,yang2014unified}, and can hardly identify community boundaries in a real network with extensive overlapping communities and hierarchical structures~\cite{cai2010graph}, leading to significant performance degeneration~\cite{gillis2013fast,luo2020highly}.
	
	In contrast, \textbf{Nonnegative Matrix Factorization (NMF)} has emerged as an effective tool for community detection due to its inherent interpretability and clustering capabilities~\cite{dao2020community,ding2006orthogonal,li2024contrastive,wang2023non}. NMF decomposes a nonnegative data matrix into two nonnegative matrices, providing a parts-based representation. Given a nonnegative matrix \( X \in \mathbb{R}_+^{n \times m} \), NMF seeks to approximate it as \( X = UV^T \), where \( U \) and \( V \) are nonnegative matrices that can be interpreted as basis and coefficient matrices. 
	
	Unlike traditional NMF, for a given $n \times n$ symmetric nonnegative similarity matrix $X$, \textbf{Symmetric Nonnegative Matrix Factorization (SymNMF)}~\cite{kuang2012symmetric} aims to find a nonnegative matrix $U \in \mathbb{R}^{n \times k}$ such that 
	\begin{align}
		\label{eq:1}
		X = UU^T, \quad \text{s.t. } U \geq 0,
	\end{align}
	where $U \geq 0$ means each element in $U$ is nonnegative, $n$ is the number of nodes in the network and $k$ is the number of communities.
	The symmetry of SymNMF naturally fits the pairwise similarity representation of network data, eliminating the need for transposition or alignment~\cite{luo2021symmetric}. Additionally, non-negativity constraints enhance interpretability, as the elements of matrix $U$ directly encode the membership degree of nodes in each community, which is useful for overlapping communities. In contrast, spectral clustering and traditional NMF often require additional post-processing steps to achieve similar functionality~\cite{xu2003document}.
	
	
	Research on community detection using SymNMF has made significant progress in recent years. The conceptual framework of SymNMF for graph clustering was established early on by Kuang et al., elucidating its equivalence and advantages relative to spectral clustering for undirected graphs~\cite{kuang2012symmetric}.Subsequent applications include rank‑2 SymNMF for hierarchical community detection in sparse networks~\cite{du2017hierarchical} and an overlapping community detection approach on ego‑splitting subgraphs~\cite{huang2021overlapping}.
	
	To enhance representation capability, researchers have incorporated graph regularization and higher‑order proximity. The SGNMF model~\cite{liu2023symmetry} uses symmetry and graph regularizers to preserve global and local network geometry. HSGN~\cite{liu2022high} integrates high‑order proximity measures. Variants such as Relaxed Symmetric NMF and Constraint‑Fusion‑induced SymNMF further balance symmetry and representation learning~\cite{liu2024relaxed}.
	
	Despite these advances, most existing SymNMF methods assume fully observed, noise-free adjacency matrices, while real-world data are often incomplete due to sampling, noise, or privacy protection~\cite{candes2010matrix,xu2003document}. In this paper, we characterize such scenarios as \textbf{mask}. Under such conditions, directly applying SymNMF would introduce estimation bias, and heuristic imputation compromises graph sparsity~\cite{candes2012exact}, leading to inaccurate community partitions~\cite{zhang2007binary}. 
	
	Inspired by matrix completion theory~\cite{recht2011simpler}, we propose a \textbf{Masked Symmetric Nonnegative Matrix Factorization (Masked SymNMF)} model motivated by the practical need for community detection in incomplete networks. Let $E \subseteq [n] \times [n]$ denote the index set corresponding to the observed entries of the symmetric adjacency matrix $X \in \mathbb{R}^{n \times n}$. We define the masking operator 
	$\mathcal{P}_E: \mathbb{R}^{n \times n} \to \mathbb{R}^{n \times n}$ which retains the values of $X$ at the observed positions in $E$ and sets all other entries to zero. Formally, the optimization objective of the SymNMF model is formulated as:
	\begin{equation}
		\label{eq:2}
		\min_{U} \frac{1}{2} \left\| \mathcal{P}_E(X - UU^T) \right\|_F^2, \quad \text{s.t. } U \geq 0,
	\end{equation}
	where the masking operator is defined by $		\left[\mathcal{P}_{E}(X) \right]_{ij} = \begin{cases}
			X_{ij}, & (i, j) \in E \\
			0, & \text{otherwise,}
		\end{cases}$
	$U\in \mathbb{R}^{n \times r}$ is a nonnegative matrix of a factorization. The proposed Masked SymNMF directly decomposes incomplete data via a masking operator, avoiding imputation bias. A regularization term further links symmetric and asymmetric decompositions, enabling tractable alternating optimization.
	
	The core problem addressed in this study is: how to design SymNMF models and efficient solution algorithms suitable for masked scenarios, and quantify the impact of key parameters on performance to better align with
	real-world scenarios in social network analysis and community detection. Solving this problem not only breaks through the application limitations of traditional methods but also provides theoretical support and practical tools for incomplete network analysis.
	

	\subsection{Contributions}
	
	\begin{itemize}
		\item  We introduce a masking operator to handle missing observations, avoiding imputation bias and improving community detection accuracy. 
		\item We reformulate the symmetric problem as an asymmetric one with a regularization term and provide a proof of the properties of the exact penalty function (\hyperlink{Thm1}{Theorem 1}), demonstrating the equivalence of the two models. 
		\item We present the Alternating Nonnegative Least Squares (ANLS) algorithm framework for the masked model with Multiplicative Update (MU), Hierarchical Alternating Least Squares (HALS), and Projected Gradient Descent (PGD) updates, prove convergence and demonstrate extensibility via additional regularization terms.
		\item We validate key parameters on synthetic data and show Masked SymNMF outperforms SymNMF across missing rates in real network.
	\end{itemize}
	
	\subsection{Outline of the Paper}
	The remainder of this paper is structured as follows: \Cref{Section2} introduces the Masked SymNMF model and presents the ANLS algorithm framework. \Cref{Section3} provides proofs of model equivalence (properties of the exact penalty function) and convergence analysis of the ANLS algorithm. \Cref{Section4} derives the update formulas for the HALS, MU, and PGD algorithms based on ANLS for both the basic symmetric and asymmetric models, as well as models with additional regularization terms, demonstrating the generalizability of the Masked SymNMF model. \Cref{Section5} presents numerical experimental results for both synthetic data and real networks. \Cref{Section6} summarizes the paper and identifies future research directions.
	
	\section{Masked SymNMF Model}
    \label{Section2}
	\subsection{From NMF to SymNMF}
	In community detection tasks, the network topology is typically represented by a symmetric adjacency matrix $X \in \mathbb{R}^{n \times n}$.
	A direct modeling approach for handling partially observed data is Masked SymNMF, as shown in Equation \eqref{eq:2}.
	Here $U \in \mathbb{R}^{n \times k}$ is the community membership matrix whose
	$i$-th row represents the nonnegative membership of node $i$ across $k$ communities. 
	
	However, Equation~\eqref{eq:2} is nonconvex in $U$ due to the \(UU^T\) coupling, making optimization challenging. First-order methods such as projected gradient descent converge slowly and are sensitive to initialization, and alternating-based algorithms (such as HALS) for nonsymmetric NMF utilize the splitting property of \eqref{eq:3} and thus cannot be used for \eqref{eq:2}~\cite{zhu2018dropping}. 
	To overcome these difficulties, we split the bilinear form of $U$ into two different factors, transforming the symmetric NMF into a nonsymmetric one and construct the following regularized asymmetric objective function:
	\begin{align}
		\label{eq:3}
		\min_{U, V} f(U, V) = \frac{1}{2} \left\| \mathcal{P}_E(X - UV^T) \right\|_F^2 + \frac{\lambda}{2} \|U - V\|_F^2, \quad \text{s.t. } U \geq 0, V \geq 0,
	\end{align}
	where $\lambda > 0$ is a regularization parameter.

    This reformulation offers two key advantages. First, it enables alternating optimization with tractable subproblems: fixing one variable in $(U,V)$ yields a convex quadratic subproblem with non-negativity constraints, for which efficient algorithms exist (e.g., projected gradient method~\cite{lin2007projected,huang2015quadratic}, alternating least squares method~\cite{kim2008nonnegative}). This makes the alternating solution process far more manageable than that of the original symmetric model $X = UU^T$. Second, the exact penalty property provides theoretical guarantees. The term $\frac{\lambda}{2} \|U - V\|_F^2$ penalizes discrepancy between $U$ and $V$; driving them to be identical at optimum recovers the original symmetric formulation. Although smooth penalties are generally not exact, we prove in \Cref{Section3} that for the masked model \eqref{eq:3} the exact penalty property remains valid. Thus the transformed model accurately recovers the solution of the original symmetric problem.

	\subsection{ANLS Algorithm Framework}
	We adopt a two-block ANLS algorithm framework, which constitutes a classical implementation of the Block Coordinate Descent (BCD) method~\cite{beck2013convergence}. Its core idea is to fix one block and solve the resulting nonnegative least-squares subproblem with respect to the other block, achieving incremental optimization of the objective function through cyclical alternating updates.
	
	The ANLS algorithm exactly solves the $U$ and $V$ subproblems, ensuring the objective decreases monotonically. These convex quadratic subproblems with nonnegativity constraints can be efficiently solved via block principal pivoting or active-set methods~\cite{kim2008nonnegative}. Our ANLS framework is general, requiring no equivalence between $U$ and $V$, and alternates over the two blocks (\hyperref[alg:2.1]{Algorithm 2.1}).
	\begin{algorithm}
		\caption{Masked SymANLS}
        \label{alg:2.1}
		\begin{algorithmic}[1]
			\State  
			\textbf{Initialization:} $\mathcal{P}_E(X)$, rank $r$, $\lambda > 0$, $U_0, V_0$
			\While{stop criterion not met}
			\State $U_k = \arg\min_{U \geq 0} \frac{1}{2} \|\mathcal{P}_E(X - UV_{k-1}^T)\|_F^2 + \frac{\lambda}{2}\|U - V_{k-1}\|_F^2$
			\State $V_k = \arg\min_{V \geq 0} \frac{1}{2} \|\mathcal{P}_E(X - U_k V^T)\|_F^2 + \frac{\lambda}{2} \|U_k - V\|_F^2$
			\State $k = k + 1$
			\EndWhile
			\State \textbf{Output:} factorization $(U_k, V_k)$
		\end{algorithmic}
	\end{algorithm}
	Within this framework, each subproblem can be precisely solved using multiple optimization algorithms, such as MU, HALS, or PGD (detailed in \Cref{Section4}). 
	
	\section{Equivalence Proof and Convergence Analysis}\label{Section3}
	\subsection{Problem Statement}
	We first recall the general nonnegative matrix factorization problem with a mask: given a similarity matrix $X \in \mathbb{R}^{n \times m}$, a factorization rank $r$, solve
	\begin{align}
		\label{eq:4}
		\min_{U, V} \frac{1}{2} \left\| \mathcal{P}_E(X - UV^T) \right\|_F^2, \quad \text{s.t. } U \geq 0, V \geq 0,
	\end{align}
	where $\mathcal{P}_E$ retains only the observed entries.
	
	When the two factors $U\in \mathbb{R}^{n \times r}$ and $V\in \mathbb{R}^{m \times r}$ are required to be identical, \eqref{eq:4} reduces to the symmetric NMF problem: $X \in \mathbb{R}^{n \times n}$
	\begin{equation}
		\label{eq:5}
		\min_U \frac{1}{2} \left\| \mathcal{P}_E(X - UU^T) \right\|_F^2, \quad \text{s.t. } U \geq 0.
	\end{equation}
	We refer to the former as the nonsymmetric case.
	
	\subsection{Problem Reformulation: Dropping Symmetry}
	We transfer the symmetric NMF problem in \eqref{eq:5} to the regularized asymmetric model
	\begin{align}
		\label{eq:6}
		\min_{U, V} f(U, V) = \frac{1}{2} \left\| \mathcal{P}_E(X - UV^T) \right\|_F^2 + \frac{\lambda}{2} \|U - V\|_F^2, \quad \text{s.t. } U \geq 0, V \geq 0.
	\end{align}
	In this section, we provide a formal guarantee (\hyperlink{Thm1}{Theorem 1}) to assure that solving \eqref{eq:6} to a critical point indeed gives a critical point of \eqref{eq:5}.
	
	\hypertarget{Lem1}{\textbf{Lemma 1.~\cite{zhu2018dropping}}} 
	For any symmetric $A \in \mathbb{R}^{n \times n}$ and PSD matrix $B \in \mathbb{R}^{n \times n}$, we have
	\begin{align*}
	\sigma_n(A) \operatorname{Tr}(B) \leq \operatorname{Tr}(AB) \leq \sigma_1(A) \operatorname{Tr}(B),
	\end{align*}
	where $\sigma_i(A)$ is the $i$-th largest eigenvalue of $A$.
	
	\hypertarget{Lem2}{\textbf{Lemma 2.}} 
	For any symmetric $A \in \mathbb{R}^{n \times n}$ and
	masking operator $\mathcal{P}_E \in \mathbb{R}^{n \times n}$, we have
	\begin{align*}
	\sigma_1(\mathcal{P}_E(A)) \leq \left\| \mathcal{P}_E(A) \right\|_F \leq \left\| A \right\|_F.
	\end{align*}
    
	\textbf{Proof.}
	(i) For any matrix $M$, we have $\|M\|_F^2 = \sum_i \sigma_i^2(M)$. Since $\sigma_i \ge 0$, it follows that $\sigma_1^2(M) \le \sum_i \sigma_i^2(M) = \|M\|_F^2$. Taking the square root and letting $M = \mathcal{P}_E(A)$, we obtain $\sigma_1(\mathcal{P}_E(A)) \le \|\mathcal{P}_E(A)\|_F$.
	
	(ii) By the definition of the masking operator, $\|\mathcal{P}_E(A)\|_F^2 = \sum_{(i,j) \in E} A_{ij}^2$. Since $E$ is a subset of all possible indices, $\sum_{(i,j) \in E} A_{ij}^2 \le \sum_{i,j=1}^n A_{ij}^2 = \|A\|_F^2$. Following the work of Candès and Recht~\cite{candes2009exact} on the contractive property of the sampling operator, it is straightforward to see that $\|\mathcal{P}_E(A)\|_F \le \|A\|_F$.

    \hypertarget{Thm1}{\textbf{Theorem 1.}} 
	Suppose that $(U^*, V^*)$ is a critical point of \eqref{eq:6} satisfying $\left\| \mathcal{P}_E(U^* V^{* T}) \right\|_F < 2\lambda + \sigma_n(\mathcal{P}_E(X))$, where $\sigma_n(\cdot)$ denotes the $n$-th largest eigenvalue. Then $U^* = V^*$ and $U^*$ is a critical point of \eqref{eq:5}.
    
	\textbf{Proof.}
	The subdifferential of $f$ is given as follows
	\begin{align}
		\label{eq:7}
		\partial_U f(U, V) = \mathcal{P}_E(UV^T-X)V + \lambda(U-V) + \partial \delta_+(U), \\
		\label{eq:8}
		\partial_V f(U, V) = \mathcal{P}_E(UV^T-X)^T U - \lambda(U-V) + \partial \delta_+(V),
	\end{align}
	where $\partial \delta_+(U) = \{G \in \mathbb{R}^{n \times r} : G \circ U = 0, G \leq 0\}$ when $U \geq 0$ and otherwise $\partial \delta_+(U) = \emptyset$. Since $(U^*, V^*)$ is a critical point of \eqref{eq:3}, it satisfies
	\begin{align}
		\label{eq:9}
		\mathcal{P}_E(U^*V^{*T}-X)V^* + \lambda(U^*-V^*)+G = 0, \\
		\label{eq:10}
		\mathcal{P}_E(U^*V^{*T}-X)^T U^* - \lambda(U^*-V^*) + H = 0,
	\end{align}
	where $G \in \partial \delta_+(U^*)$ and $H \in \partial \delta_+(V^*)$. Subtracting the second from the first, we have
	\begin{align}
		\label{eq:11}
		(2\lambda I + \mathcal{P}_E(X))(U^* - V^*) = \mathcal{P}_E(V^* U^{*T})U^* - \mathcal{P}_E(U^* V^{*T})V^* - G + H,
	\end{align}
	where we utilize the fact that $X$ is symmetric, i.e., $X = X^T$. Taking the inner product of $U^* - V^*$ with both sides of the above equation gives
	\begin{equation}
		\label{eq:12}
        \begin{aligned}
		&\quad (2\lambda I + \mathcal{P}_E(X)),\, (U^* - V^*)(U^* - V^*)^T \rangle \\
		& = \mathcal{P}_E(V^* U^{*T})U^* - \mathcal{P}_E(U^* V^{*T})V^* - G + H,\, U^* - V^* \rangle.
        \end{aligned}
	\end{equation}
	In what follows, by choosing sufficiently large $\lambda$, we show that $(U^*, V^*)$ satisfying \eqref{eq:9}--\eqref{eq:10} must satisfy $U^* = V^*$. To that end, we first provide the lower bound and the upper bound for the LHS and RHS of \eqref{eq:10}, respectively. Specifically,
	\begin{equation}
		\label{eq:13}
		\begin{aligned}
		  &\quad (2\lambda I + \mathcal{P}_E(X)), (U^*-V^*)(U^*-V^*)^T \rangle \\
		& \geq \sigma_n(2\lambda I + \mathcal{P}_E(X)) \|U^*-V^*\|_F^2 \\
		& = (2\lambda + \sigma_n(\mathcal{P}_E(X))) \|U^*-V^*\|_F^2,
		\end{aligned}
	\end{equation}
	where the inequality follows from \hyperlink{Lem1}{Lemma 1}. On the other hand,
	\begin{align*}
		&\quad \mathcal{P}_E(V^*U^{*T})U^* - \mathcal{P}_E(U^*V^{*T})V^* - G + H,\ U^* - V^* \rangle \\
		&\leq \langle \mathcal{P}_E(V^*U^{*T})U^* - \mathcal{P}_E(U^*V^{*T})V^*,\ U^* - V^* \rangle \\
		&= \left\langle \tfrac{1}{2}\!\left[\mathcal{P}_E(V^*U^{*T}) + \mathcal{P}_E(U^*V^{*T})\right],\ (U^* - V^*)(U^* - V^*)^T \right\rangle \\
		&\quad - \tfrac{1}{2} \left\langle \mathcal{P}_E(U^*V^{*T}) - \mathcal{P}_E(V^*U^{*T}),\ U^*V^{*T} - V^*U^{*T} \right\rangle \\
		&\leq \left\langle \tfrac{1}{2}\!\left[\mathcal{P}_E(V^*U^{*T}) + \mathcal{P}_E(U^*V^{*T})\right],\ (U^* - V^*)(U^* - V^*)^T \right\rangle \\
		&\leq \sigma_1\!\left\{ \tfrac{1}{2}\!\left[\mathcal{P}_E(V^*U^{*T}) + \mathcal{P}_E(U^*V^{*T})\right] \right\} \left\| U^* - V^* \right\|_F^2,
	\end{align*}
	where the last inequality utilizes \hyperlink{Lem1}{Lemma 1} and the first inequality follows because $V^*, U^* \geq 0$ indicating that
	\begin{equation}
		\label{eq:14}
		-\langle G, U^* - V^* \rangle \leq 0, \quad \langle H, U^* - V^* \rangle \leq 0.
	\end{equation}
	Substituting \eqref{eq:13} and \eqref{eq:14} into \eqref{eq:12} and utilizing the assumption that $\left\| \mathcal{P}_E(U^*V^{*T}) \right\|_F \leq \alpha$, we have
    \begin{equation}
    \label{eq:15}
    \begin{aligned}
    &\quad (2\lambda + \sigma_n(\mathcal{P}_E(X))) \|U^* - V^*\|_F^2 \\
    & \leq \sigma_1 \left\{ \frac{1}{2}[\mathcal{P}_E(V^*U^{*T}) + \mathcal{P}_E(U^*V^{*T})] \right\} \|U^* - V^*\|_F^2 \\
    & \leq \alpha \|U^* - V^*\|_F^2,
    \end{aligned}
    \end{equation}
	where the last inequality utilizes \hyperlink{Lem2}{Lemma 2}. This implies that if we choose $2\lambda > \alpha - \sigma_n(\mathcal{P}_E(X))$, then $U^* = V^*$ must hold. Substituting this identity into \eqref{eq:7} gives
	\begin{equation}
		\label{eq:16}
		0 \in \mathcal{P}_E(U^*(U^*)^T - X)U^* + \partial \delta_+(U^*),
	\end{equation}
	which implies $U^*$ is a critical point of \eqref{eq:5}.
	
	\subsection{Convergence Analysis of ANLS Algorithm}
	In this section, we prove that the ANLS algorithm converges to a critical point (stationary point) of model \eqref{eq:3}. The proof is based on sufficient descent, boundedness of the sequence, and fixed-point analysis, similar to the proof for standard NMF in~\cite{kim2008nonnegative}.
	
	\hypertarget{Lem3}{\textbf{Lemma 3 (Sufficient Decrease Property)}} 
	The iterative sequence $\{(U^k, V^k)\}$ generated by the ANLS framework (\hyperref[alg:2.1]{Algorithm 2.1}) satisfies the non-strictly decreasing property of the objective function,
	$$
	f({U}_{k+1}, {V}_{k+1}) \leq f({U}_{k+1}, {V}_k) \leq f({U}_k, {V}_k), \quad \forall k \geq 0,
	$$
	with equality if and only if $({U}_k, {V}_k)$ is a fixed point of the algorithm.
	
	\textbf{Proof.} 
	By the iteration rule of ANLS, ${U}_{k+1} = \arg\min_{{U} \geq 0} f({U}, {V}_k)$, which implies ${U}_{k+1}$ is the global optimal solution of the ${U}$-subproblem with fixed ${V}_k$. Thus
	$f({U}_{k+1}, {V}_k) \leq f({U}_k, {V}_k).$
	Equality holds if and only if ${U}_k = {U}_{k+1}$.
	Similarly, ${V}_{k+1} = \arg\min_{{V} \geq 0} f({U}_{k+1}, {V})$, so
	$f({U}_{k+1},{V}_{k+1}) \leq f({U}_{k+1}, {V}_k).$
	Equality holds if and only if ${V}_k = {V}_{k+1}$.
	Combining the above two inequalities, we get $f({U}_{k+1}, {V}_{k+1}) \leq f({U}_k, {V}_k)$, with equality if and only if ${U}_k = {U}_{k+1}$ and ${V}_k = {V}_{k+1}$ (i.e., $({U}_k, {V}_k)$ is a fixed point).

    The sufficient decrease property directly follows the proof in Beck et al~\cite{beck2013convergence}, which remains valid under the masking operator $\mathcal{P}_E$.

    \hypertarget{Lem4}{\textbf{Lemma 4 (Boundedness of Iterative Sequence)}} 
	For any local search algorithm solving \eqref{eq:3} with initialization $V_0 = U_0$, suppose it sequentially decreases the objective value and the diagonal elements of the masking operator $\mathcal{P}_E$ is non-zero. Then, for any $k \geq 0$, the iterate $(U_k, V_k)$ generated by this algorithm satisfies
	\begin{equation}
		\label{eq:17}
		\begin{aligned}
			\left\| U_k \right\|_F^2 + \left\| V_k \right\|_F^2 &\leq \left( \frac{1}{\lambda} + 2\sqrt{n} \right) \left\| \mathcal{P}_E(X - U_0 U_0^T) \right\|_F^2 + 2\sqrt{n} \|\mathcal{P}_E(X)\|_F := B_0, \\
			\left\| \mathcal{P}_E(U_k V_k^T) \right\|_F &\leq \left\|\mathcal{P}_E(X - U_0 V_0^T) \right\|_F + \|\mathcal{P}_E(X)\|_F.
		\end{aligned}
	\end{equation}
	
	\textbf{Proof.}
	By the assumption that the algorithm decreases the objective function, we have
	$\frac{1}{2} \left\| \mathcal{P}_E(X - U_k V_k^T) \right\|_F^2 + \frac{\lambda}{2} \left\| U_k - V_k \right\|_F^2 \leq \frac{1}{2} \left\| \mathcal{P}_E(X - U_0 U_0^T) \right\|_F^2,$
	which further implies that
	\begin{align*}
	\begin{cases}
		\left\| \mathcal{P}_E(X - U_k V_k^T) \right\|_F \leq \left\| \mathcal{P}_E(X - U_0 U_0^T)\right\|_F \\
		\frac{\lambda}{2} \left( \left\| U_k \right\|_F^2 + \left\| V_k \right\|_F^2 - 2 \left| \left\langle U_k^T V_k, I_r \right\rangle \right| \right) \leq \frac{\lambda}{2} \left\| U_k - V_k \right\|_F^2 \leq \frac{1}{2} \left\| \mathcal{P}_E(X - U_0 U_0^T) \right\|_F^2.
	\end{cases}
	\end{align*}
     From the first line, we obtain$\left\| \mathcal{P}_E(U_k V_k^T) \right\|_F \leq \left\|\mathcal{P}_E(X - U_0 V_0^T)\right\|_F + \|\mathcal{P}_E(X)\|_F,$
	while the second line leads to $\left\| U_k \right\|_F^2 + \left\| V_k \right\|_F^2\leq \frac{1}{\lambda} \left\| \mathcal{P}_E(X - U_0 U_0^T) \right\|_F^2 + 2\mathrm{Tr}(U_k V_k^T).$
    
    Based on our assumption that the diagonal elements of the masking operator $\mathcal{P}_E$ are non-zero, the set of observed entries $E$ contains all diagonal indices $(i,i)$ for $i=1,\dots,n$. Since $U_k \ge 0$ and $V_k \ge 0$, the diagonal elements $(U_k V_k^T)_{ii}$ are non-negative. By applying the Cauchy-Schwarz inequality, we can bound the trace term using the masked Frobenius norm:
    $\text{Tr}(U_k V_k^T)\le \sqrt{n} \sqrt{\sum_{i=1}^n (U_k V_k^T)_{ii}^2}\le \sqrt{n} \sqrt{\sum_{(i,j)\in E} (U_k V_k^T)_{ij}^2}.$
    
    Substituting this inequality and the upper bound of $||\mathcal{P}_E(U_k V_k^T)||_F$ into the equation, we obtain a rigorous upper bound $B_0$:
    $||U_k||_F^2 + ||V_k||_F^2\le \frac{1}{\lambda}||\mathcal{P}_E(X - U_0 U_0^T)||_F^2 + 2\sqrt{n} \left( ||\mathcal{P}_E(X - U_0 V_0^T)||_F + ||\mathcal{P}_E(X)||_F \right) := B_0.$
    
	
	From \eqref{eq:17}, two facts about the iteration can be derived: The first equation in \eqref{eq:17} implies that both $U_k$ and $V_k$ are bounded, with upper bounds decay as $\lambda$ increases.
	The second equality in \eqref{eq:17} indicates that it is upper bounded by a quantity independent of $\lambda$. This reveals a result that if the iterative algorithm converges and the iterates $(U_k, V_k)$ converge to a critical point $(U^*, V^*)$, then $U^* V^{*T}$ is also bounded, independent of the value of $\lambda$. Together with \hyperlink{Thm1}{Theorem 1}, this guarantees that by choosing sufficiently large $\lambda$, many local search algorithms can be utilized to find critical points of \eqref{eq:5}.
    
    The boundedness analysis in \hyperlink{Lem4}{Lemma 4} follows the framework established by Zhu et al.~\cite{zhu2018dropping} for SymNMF, extended here to accommodate the masking operator $\mathcal{P}_E$.

    \hypertarget{Lem5}{\textbf{Lemma 5 (Convex Quadratic Subproblems)}}
	The subproblems of solving 
	${U}$ with fixed ${V}$ or solving ${V}$ with fixed ${U}$ are convex quadratic optimization problems with non-negativity constraints, and there exists a unique global optimal solution.
     
    \textbf{Proof.}
    Consider the subproblem for $U$ (with $V$ fixed):
    \begin{align*}
    \min_U \frac{1}{2} \bigl\| \mathcal{P}_E(X - U V^T) \bigr\|_F^2 + \frac{\lambda}{2} \| U - V \|_F^2,\quad \text{s.t. } U \geq 0.\end{align*} Vectorizing $U$, the objective is a quadratic form $  \frac{1}{2} \operatorname{vec}(U)^T H \operatorname{vec}(U) + b^T \operatorname{vec}(U) + c  $, where the Hessian $  H  $ is block-structured as $  H = \sum_{(i,j) \in E} (v_j v_j^T) \otimes e_i e_i^T + \lambda I $ (with $e_i$ the standard basis vector). Each term $ (v_j v_j^T) \otimes e_i e_i^T $ is positive semidefinite, and the regularization term $  \lambda I  $ ($ \lambda > 0 $) makes $H$ positive definite. Hence the objective is strictly convex over the convex feasible set $  U \geq 0  $, guaranteeing a unique global minimizer. The subproblem for $V$ (with $U$ fixed) follows symmetrically.

    The convexity and uniqueness of ANLS subproblems are well-established. Kim and Park~\cite{kim2008nonnegative} proved their convexity and unique solvability; Lin~\cite{lin2007projected} and Huang et al.~\cite{huang2015quadratic} analyzed related projected gradient methods. The masking operator \(\mathcal{P}_E\) only restricts summation to observed entries and preserves Hessian definiteness, so the above proof follows directly.
    
    \hypertarget{Thm2}{\textbf{Theorem 2.}}
	Choose $\lambda > \frac{1}{2}\bigl(\|{\mathcal{P}_E(X)}\|_2 + \bigl\|\mathcal{P}_E({X} - {U}_0 {U}_0^T)\|_F - \sigma_n({\mathcal{P}_E(X}))\bigr)$ for \eqref{eq:6}, where $\|{X}\|_2$ means the largest singular value of $X$. For any local search algorithm solving \eqref{eq:6} with initialization $V_0 = U_0$, if it sequentially decreases the objective value, is convergent, and converges to a critical point $({U}^*, {V}^*)$ of \eqref{eq:6}, then we have ${U}^* = {V}^*$ and that ${U}^*$ is also a critical point of \eqref{eq:5}.

    \textbf{Proof.} 
    By \hyperlink{Lem4}{Lemma 4}, the sequence ${(U_k,V_k)}$ generated by any algorithm is bounded, and $\left\| \mathcal{P}_E(U_k V_k^T) \right\|_F$ is upper bounded by a constant independent of $\lambda$. Let $(U^*, V^*)$ be the limit point, then we have
    $\left\| \mathcal{P}_E(U^* V^{* T}) \right\|_F \leq \left\| \mathcal{P}_E(X-U_0 V_0^T) \right\|_F + \left\| \mathcal{P}_E(X) \right\|_F.$
    The chosen lower bound on $\lambda$ ensures that the condition $\left\| \mathcal{P}_E(U^* V^{* T}) \right\|_F < 2\lambda + \sigma_n(\mathcal{P}_E(X))$ in \hyperlink{Thm1}{Theorem 1} is satisfied. Since $(U^*, V^*)$ is a critical point of \eqref{eq:6} and satisfies this upper bound, it directly follows from \hyperlink{Thm1}{Theorem 1} that $U^* = V^*$ and that $U^*$ is a critical point of the original symmetric problem \eqref{eq:5}.
    
	\hyperlink{Thm2}{Theorem 2} indicates that instead of directly solving \eqref{eq:5}, we can solve \eqref{eq:6} with a sufficiently large regularization parameter $\lambda$. The latter is similar to the nonsymmetric NMF \eqref{eq:5} and exhibits similar splitting properties, enabling us to utilize efficient alternating-type algorithms. In the next section, we propose an alternating-based algorithm for handling \eqref{eq:6}.
    
	\hypertarget{Thm3}{\textbf{Theorem 3(Convergence of \hyperref[alg:2.1]{Algorithm 2.1}).}} Let $\{(U^k, V^k)\}$ be the sequence generated by the ANLS framework. Then the objective values decrease monotonically and are bounded below, the sequence itself is bounded (by \hyperlink{Lem3}{Lemma 3}), and every limit point $(U^*,V^*)$ satisfies the KKT conditions of model \eqref{eq:6}, hence is a stationary point.
	
	\textbf{Proof.}
	According to the sufficient decrease property shown in \hyperlink{Lem3}{Lemma 3}, the ANLS framework ensures $f(U_{k+1}, V_{k+1}) \le f(U_{k+1}, V_k) \le f(U_k, V_k)$. Thus, the sequence $\{f(U_k, V_k)\}$ is monotonically decreasing. Since the objective function is bounded below by 0, the sequence converges.
	
	Following \hyperlink{Lem4}{Lemma 4}, the objective function's descent ensures that the sequence $\{(U_k, V_k)\}$ remains within a compact set defined by $\|U\|_F^2 + \|V\|_F^2 \le B_0$. By the Bolzano-Weierstrass theorem, there exists at least one limit point $(U^*, V^*)$.
	
	Since the subproblems are convex quadratic (\hyperlink{Lem5}{Lemma 5}) and the objective function is continuously differentiable, we apply the standard two-block coordinate descent theory. At the limit point $(U^*, V^*)$, the objective function cannot be further decreased by updating either $U$ or $V$. Consequently, the partial derivatives satisfy the KKT optimality conditions: $\nabla_U f(U^*, V^*) \ge 0$ and $\nabla_V f(U^*, V^*) \ge 0$, along with complementary slackness.
	Thus, $(U^*, V^*)$ is a stationary point of the problem.
	
	\hyperlink{Thm3}{Theorem 3} demonstrates that the Masked SymNMF model with masking operator $\mathcal{P}_E$ satisfies sufficient descent and boundedness of the iteration sequence for its ANLS algorithm. The iteration sequence generated by the algorithm must converge to a critical point of the objective function.
	In fact, the masking operator $\mathcal{P}_E$ constrains only the observed entries without altering the smoothness of the objective function or the convex quadratic nature of the subproblems. Consequently, the convergence theory of the classical two-block BCD can be directly extended to the masked model.

	\section{Modified Algorithms for Masked SymNMF}\label{Section4}
	In this section, we first present the Masked SymNMF models with and without common regularization terms. For both the original symmetric model and the transformed asymmetric model, we derive the iterative update rules of three classical algorithms: MU~\cite{lee1999learning}, HALS~\cite{kim2008nonnegative,zhu2018dropping}, and PGD~\cite{cichocki2009fast,huang2015quadratic} based on ANLS. These algorithms have been widely validated in standard NMF~\cite{cichocki2007hierarchical,lee2000algorithms}, and we adapt them to the Masked SymNMF framework here. On this basis, we further introduce graph regularization widely used in community detection~\cite{cai2010graph,liu2023symmetry,yang2014unified}, as well as general $\ell_1$ regularization~\cite{hoyer2004non,tibshirani1996regression} and orthogonality constraints~\cite{choi2008algorithms,ding2005equivalence}, to construct enhanced regularized models and present their corresponding algorithmic updates. This provides a complete and extensible algorithmic framework for network community detection under masking.
	
	\subsection{Basic Models}
	In this section, we present the objective functions of the original masked symmetric model ($O_1$) and the transformed asymmetric regularized model ($O_2$), then derive the update rules of three core solving algorithms based on the ANLS framework. MU rule is widely adopted in NMF problems due to its simple implementation and automatic non‑negativity preservation~\cite{lee1999learning}. The HALS accelerates convergence via column‑wise optimization~\cite{cichocki2009fast,cichocki2007hierarchical}. The PGD method ensures iterative stability using gradient information~\cite{huang2015quadratic,kim2008nonnegative}. All the above algorithms adapt to the masking operator $\mathcal{P}_E$ and can directly handle incomplete network adjacency matrices.
	
	As mentioned earlier, the objective functions corresponding to NMF and SymNMF with mask are $O_1$ and $O_2$. For the objective function $O_1$, we adopt the MU rule for optimization, and its update formula is given as \eqref{eq:19}. 
	\begin{align}
		\label{eq:18}
		O_1 = \frac{1}{2} \left\| \mathcal{P}_E(X - UU^T) \right\|_F^2, \quad \text{s.t.} \ U \geq 0.
	\end{align}
	\begin{align}
    \textbf{[MU]}:
		\label{eq:19}
		u_{ik} &\leftarrow u_{ik} \frac{[\mathcal{P}_E(X)U]_{ik}}{[\mathcal{P}_E(UU^T)U]_{ik}}.
	\end{align}
	
	For the asymmetric model $O_2$, the MU rules act on the factor matrices $U$ and $V$ separately, as shown in Equations \eqref{eq:21} and \eqref{eq:22}.
	\begin{align}
		\label{eq:20}
		O_2 = \frac{1}{2} \left\| \mathcal{P}_E(X - UV^T) \right\|_F^2 + \frac{\lambda}{2} \|U - V\|_F^2, \quad \text{s.t.} \ U \geq 0, V \geq 0.
	\end{align}
	\begin{align}
    \textbf{[MU]}:
		\label{eq:21}
		u_{ik} &\leftarrow u_{ik} \frac{[\mathcal{P}_E(X)V + \lambda V]_{ik}}{[\mathcal{P}_E(UV^T)V + \lambda U]_{ik}}, \\
		\label{eq:22}
		v_{jk} &\leftarrow v_{jk} \frac{[\mathcal{P}_E(X^T)U + \lambda U]_{jk}}{[\mathcal{P}_E(VU^T)U + \lambda V]_{jk}}.
	\end{align}
    
	In addition to MU, the HALS algorithm further accelerates convergence by updating columns sequentially, yielding the following update rules.
    \begin{align}
    \textbf{[HALS]}:
    \label{eq:23}
    u_{ik} &\leftarrow u_{ik} \frac{[\mathcal{P}_E(X_k){v}_k + \lambda{v}_k]_i}{[\mathcal{P}_E({v}_k {v}_k^T)\mathbf{1} + \lambda \mathbf{1}]_i}, \\
     \label{eq:24}
    v_{jk} &\leftarrow v_{jk} \frac{[\mathcal{P}_E({X}_k^T) {u}_k + \lambda {u}_k]_j}{[\mathcal{P}_E({u}_k {u}_k^T) \mathbf{1} + \lambda \mathbf{1}]_j},
    \end{align}
	where $X_k = X - \sum_{j \ne k} u_j v_j^T$.
	
	Finally, the gradient calculation formulas of PGD based on the ANLS framework are provided, offering a core basis for solving with gradient‑based algorithms.
	\begin{align}
    \textbf{[PGD]}:
		\label{eq:25}
		\nabla_{\!U}f = (\mathcal{P}_E(UV^T-X))V + \lambda(U-V), \\
		\label{eq:26}
		\nabla_{\!V}f = (\mathcal{P}_E(VU^T-X))U + \lambda(V-U).
	\end{align}
	All updates satisfy non‑negativity constraints and are compatible with the masked objective (see the equivalence proof in \Cref{Section3}).
	
	\subsection{Models with Regularization Terms}
	
	To enhance performance on sparse or noisy networks, this section introduces three types of regularization constraints on the basic model: graph regularization preserves the local topological structure of the network and is a mainstream enhancement for community detection~\cite{cai2010graph,liu2023symmetry}; $\ell_1$ regularization induces sparsity in the membership matrix to improve the interpretability of community assignments~\cite{hoyer2004non,tibshirani1996regression}; orthogonality constraints enhance the discrimination between communities and alleviate ambiguous overlapping~\cite{choi2008algorithms,ding2005equivalence}. 
	For each case, we provide the objective and corresponding update rules for both symmetric ($O_{11}$, $O_{12}$, $O_{13}$) and asymmetric ($O_{21}$, $O_{22}$, $O_{23}$) models. 
	
	\subsubsection{Graph Regularization}
	Graph regularization preserves the local geometric structure (manifold structure) of the network and improves community detection accuracy, especially for sparse or noisy social networks~\cite{cai2010graph}. In this subsection, we present the graph‑regularized symmetric model $O_{11}$ and asymmetric model $O_{21}$, as well as the corresponding update formulas.
	
	The objective function of the symmetric model $O_{11}$ with graph regularization is given, integrating the graph regularization term into the symmetric factorization framework. The MU rule is derived based on this model, as shown in Equation \eqref{eq:28}.
	\begin{align}
		\label{eq:27}
		O_{11} = \frac{1}{2} \left\| \mathcal{P}_E(X - UU^T) \right\|_F^2 + \frac{\alpha}{2}\mathrm{Tr}\left(U^T L U\right), \quad \text{s.t.} \ U \geq 0.
	\end{align}
	\begin{align}
    \textbf{[MU]}:
		\label{eq:28}
		u_{ik} &\leftarrow u_{ik} \frac{[\mathcal{P}_E(X)U + \alpha WU]_{ik}}{[\mathcal{P}_E(UU^T)U + \alpha DU]_{ik}}.
	\end{align}
		
	For the asymmetric model $O_{21}$, the update formulas of MU, HALS and PGD based on ANLS are given in sequence.
	\begin{align}
		\label{eq:29}
		O_{21} = \frac{1}{2} \left\| \mathcal{P}_E(X - UV^T) \right\|_F^2 + \frac{\lambda}{2} \|U - V\|_F^2 + \frac{\alpha}{2}\mathrm{Tr}\left(V^T L V\right), \quad \text{s.t.} \ U \geq 0, V \geq 0.
	\end{align}
	
	The corresponding MU rules need to handle both the penalty term and the graph regularization, and are given as follows.
	\begin{align}
    \textbf{[MU]}:
		\label{eq:30}
		u_{ik} &\leftarrow u_{ik} \frac{[\mathcal{P}_E(X)V + \lambda V]_{ik}}{[\mathcal{P}_E(UV^T)V + \lambda U]_{ik}}, \\
		\label{eq:31}
		v_{jk} &\leftarrow v_{jk} \frac{[\mathcal{P}_E(X^T)U + \lambda U + \alpha WV]_{jk}}{[\mathcal{P}_E(VU^T)U + \lambda V + \alpha DV]_{jk}}.
	\end{align}
	
	For the same model, the column‑wise updates of the HALS algorithm can also incorporate graph regularization information.
	\begin{align}
    \textbf{[HALS]}:
		\label{eq:32}
        u_{ik} &\leftarrow u_{ik} \frac{[\mathcal{P}_E({X}_k) {v}_i + \lambda{v}_i]_k}{[\mathcal{P}_E({v}_k {v}_k^T) \mathbf{1} + \lambda \mathbf{1}]_i}, \\
		\label{eq:33}
        v_{jk} &\leftarrow v_{jk} \frac{[\mathcal{P}_E({X}_k^T) {u}_k + \lambda{u}_k + \alpha W {v}_k]_j}{[\mathcal{P}_E({u}_k {u}_k^T) \mathbf{1} + \lambda \mathbf{1} + \alpha d]_j},
	\end{align}
	where $X_k = X - \sum_{j \ne k} u_j v_j^T$, and $d$ is the degree vector (i.e. $d_j = \sum_k X_{jk}$).
	
	In the ANLS‑based PGD method, the gradient update rules are given by Equations \eqref{eq:34} and \eqref{eq:35}.
	\begin{align}
    \textbf{[PGD]}:
		\label{eq:34}
		\nabla_{\!U}F = (\mathcal{P}_E(UV^T-X))V + \lambda(U-V), \\
		\label{eq:35}
		\nabla_{\!V}F = (\mathcal{P}_E(VU^T-X))U + \lambda(V-U) + \alpha LV.
	\end{align}
	
	\subsubsection{$\ell_1$ Regularization}
	
	$\ell_1$ regularization induces sparsity in the factor matrices, which is particularly suitable for overlapping community detection where node memberships are sparse~\cite{hoyer2004non}. This subsection presents the $\ell_1$‑regularized symmetric model $O_{12}$ and asymmetric model $O_{22}$ with corresponding update rules.
	
	$\ell_1$ regularization is given to achieve factor matrix sparsity via $\ell_1$ norm constraint.
	\begin{align}
		\label{eq:36}
		O_{12} = \frac{1}{2} \left\| \mathcal{P}_E(X - UU^T) \right\|_F^2 + \beta \|U\|_1, \quad \text{s.t.} \ U \geq 0,
	\end{align}
	where $\beta > 0$ controls the sparsity strength. The MU update rule is obtained based on this model, adding the $\ell_1$ regularization term to the denominator to induce a sparse solution. The MU update rules for $O_{12}$ is
	\begin{align}
    \textbf{[MU]}:
		\label{eq:37}
		u_{ik} &\leftarrow u_{ik} \frac{[\mathcal{P}_E(X)U]_{ik}}{[\mathcal{P}_E(UU^T)U + \beta \mathbf{1}_{n \times r}]_{ik}}.
	\end{align}
	
	Incorporating sparse regularization into the asymmetric model yields $O_{22}$, where the penalty term $\frac{\lambda}{2}\|U-V\|_F^2$ and the $\ell_1$ term jointly act on $U$.
	\begin{align}
		\label{eq:38}
		O_{22} = \frac{1}{2} \left\| \mathcal{P}_E(X - UV^T) \right\|_F^2 + \frac{\lambda}{2} \|U - V\|_F^2 + \beta \|U\|_1, \quad \text{s.t.} \ U \geq 0, V \geq 0.
	\end{align}
	
	The MU update rules for $O_{22}$ are
	\begin{align}
    \textbf{[MU]}:
		\label{eq:39}
		u_{ik} &\leftarrow u_{ik} \frac{[\mathcal{P}_E(X)V + \lambda V]_{ik}}{[\mathcal{P}_E(UV^T)V + \lambda U + \beta \mathbf{1}_{n \times r}]_{ik}}, \\
		\label{eq:40}
		v_{jk} &\leftarrow v_{jk} \frac{[\mathcal{P}_E(X^T)U + \lambda U]_{jk}}{[\mathcal{P}_E(VU^T)U + \lambda V]_{jk}}.
	\end{align}
	
	The corresponding gradient expressions in the ANLS-based PGD method are
	\begin{align}
    \textbf{[PGD]}:
		\label{eq:41}
		\nabla_{\!U}F = (\mathcal{P}_E(UV^T-X))V + \lambda(U-V) + \beta \, \text{sign}(U), \\
		\label{eq:42}
		\nabla_{\!V}F = (\mathcal{P}_E(VU^T-X))U + \lambda(V-U).
	\end{align}
	
	\subsubsection{Orthogonality Constraint}
	Orthogonality constraints force community basis vectors to be approximately orthogonal, improving the separation of community partitions, which is suitable for non‑overlapping or weakly overlapping community detection~\cite{ding2006orthogonal, li2025orthogonal}. This subsection presents the orthogonality‑constrained symmetric model $O_{13}$ and asymmetric model $O_{23}$ with corresponding algorithm iterations.

	The objective function of the symmetric orthogonality-constrained model $O_{13}$ is
	\begin{align}
		\label{eq:43}
		O_{13} = \frac{1}{2} \left\| \mathcal{P}_E(X - UU^T) \right\|_F^2 + \frac{\gamma}{4}\left\|U^T U - I \right\|_F^2, \quad \text{s.t.} \ U \geq 0,
	\end{align}
	where $\gamma > 0$ is the orthogonality penalty parameter. The MU rule adds a term $\gamma U$ in the denominator:
	\begin{align}
    \textbf{[MU]}:
		\label{eq:44}
		u_{ik} &\leftarrow u_{ik} \frac{[\mathcal{P}_E(X)U+ \gamma U]_{ik}}{[\mathcal{P}_E(UU^T)U + \gamma UU^TU]_{ik}}.
	\end{align}
	
	For the asymmetric orthogonality-constrained model $O_{23}$, the orthogonality penalty is imposed on $U$, as shown in Equation \eqref{eq:45}
	\begin{align}
		\label{eq:45}
		O_{23} = \frac{1}{2} \left\| \mathcal{P}_E(X - UV^T) \right\|_F^2 + \frac{\lambda}{2} \|U - V\|_F^2 + \frac{\gamma}{4}\left\|U^T U - I \right\|_F^2, \quad \text{s.t.} \ U \geq 0, V \geq 0.
	\end{align}
	
	The MU rules for $O_{23}$ are
	\begin{align}
    \textbf{[MU]}:
		\label{eq:46}
		u_{ik} &\leftarrow u_{ik} \frac{[\mathcal{P}_E(X)V + \lambda V + \gamma U]_{ik}}{[\mathcal{P}_E(UV^T)V + \lambda U + \gamma UU^TU]_{ik}}, \\
		\label{eq:47}
		v_{jk} &\leftarrow v_{jk} \frac{[\mathcal{P}_E(X^T)U + \lambda U]_{jk}}{[\mathcal{P}_E(VU^T)U + \lambda V]_{jk}}.
	\end{align}
	
	In the ANLS-based PGD method, the orthogonality constraint contributes an extra gradient term $\gamma U$. The modified gradient expressions are
	\begin{align}
    \textbf{[PGD]}:
		\label{eq:48}
		\nabla_{\!U}F = (\mathcal{P}_E(UV^T-X))V + \lambda(U-V) + \gamma U(U^TU-I), \\
		\label{eq:49}
		\nabla_{\!V}F = (\mathcal{P}_E(VU^T-X))U + \lambda(V-U).
	\end{align}
	
	The derivations follow the same pattern as the basic models, substituting the additional penalty terms into the gradients. This yields a complete, extensible algorithmic framework for masked community detection. Subsequent numerical experiments will evaluate graph regularization for community detection.
	
	\section{Numerical Experiments}\label{Section5}
	In this section, we evaluate the performance of the proposed Masked SymNMF model and its corresponding optimization algorithms (e.g., HALS, MU, PGD based on ANLS).
	
	\subsection{Synthetic Datasets}
	To systematically examine the theoretical properties and algorithmic behavior of the proposed Masked SymNMF model, we first conduct experiments on synthetic data.
	
	\subsubsection{Data Generation}
    The synthetic matrices are generated as follows. Random nonnegative initial factors \( U_0, V_0 \in \mathbb{R}^{m \times r} \) are uniformly sampled from the interval \((0, 1)\). A symmetric mask matrix \({P}_E \) is constructed by symmetrizing a random matrix and thresholding at a given mask rate. A low-rank component \( M = AA^T \) is formed with a nonnegative random matrix \( A \). Optionally, Gaussian noise may be added. The observed matrix is obtained by elementwise application of the mask:$M \leftarrow P_E \odot M$, where $\odot$ denotes the Hadamard product.
    
	
	\subsubsection{Impact of Matrix Dimension and Density}
	
\begin{figure}[ht]
		\centering
		
		\begin{subfigure}{0.45\textwidth}
			\centering
			\includegraphics[width=\textwidth]{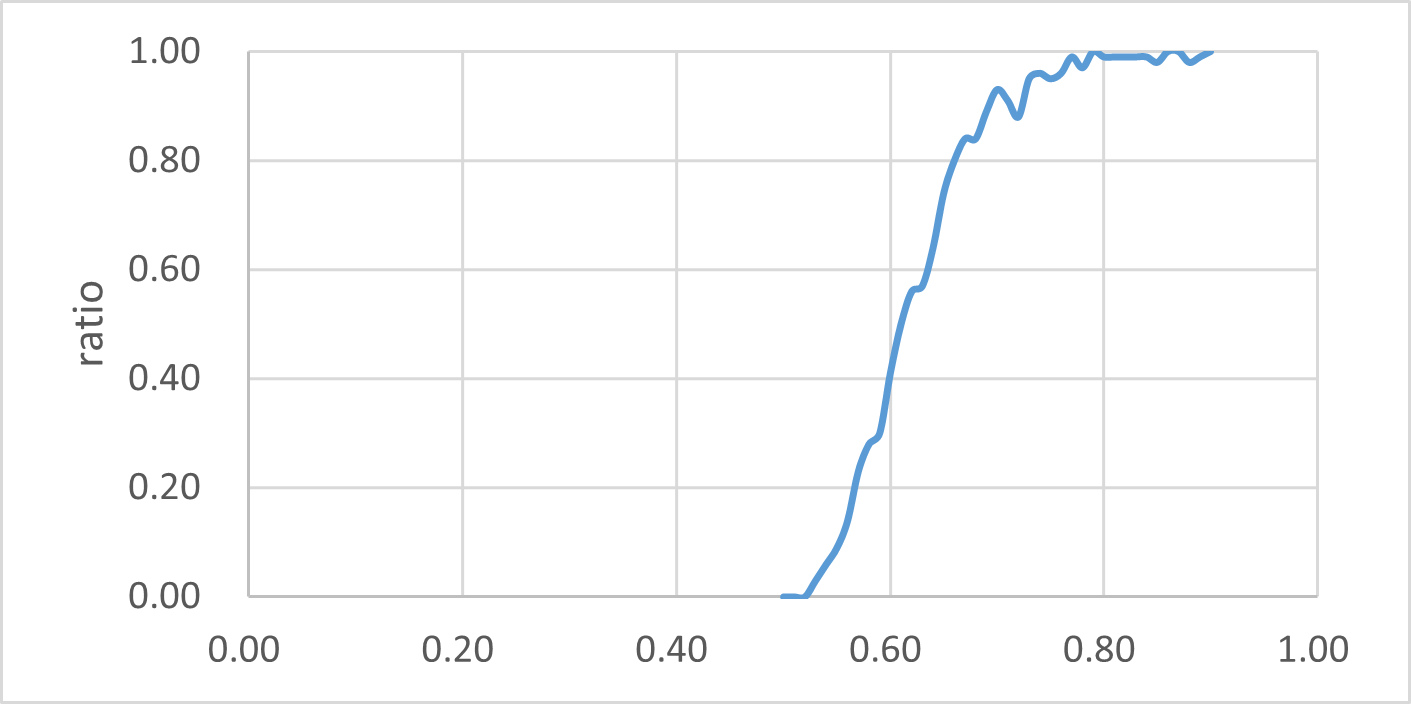}  
			\caption{m=100}
			\label{fig:(a)}
		\end{subfigure}
		\hfill
		\begin{subfigure}{0.45\textwidth}
			\centering
			\includegraphics[width=\textwidth]{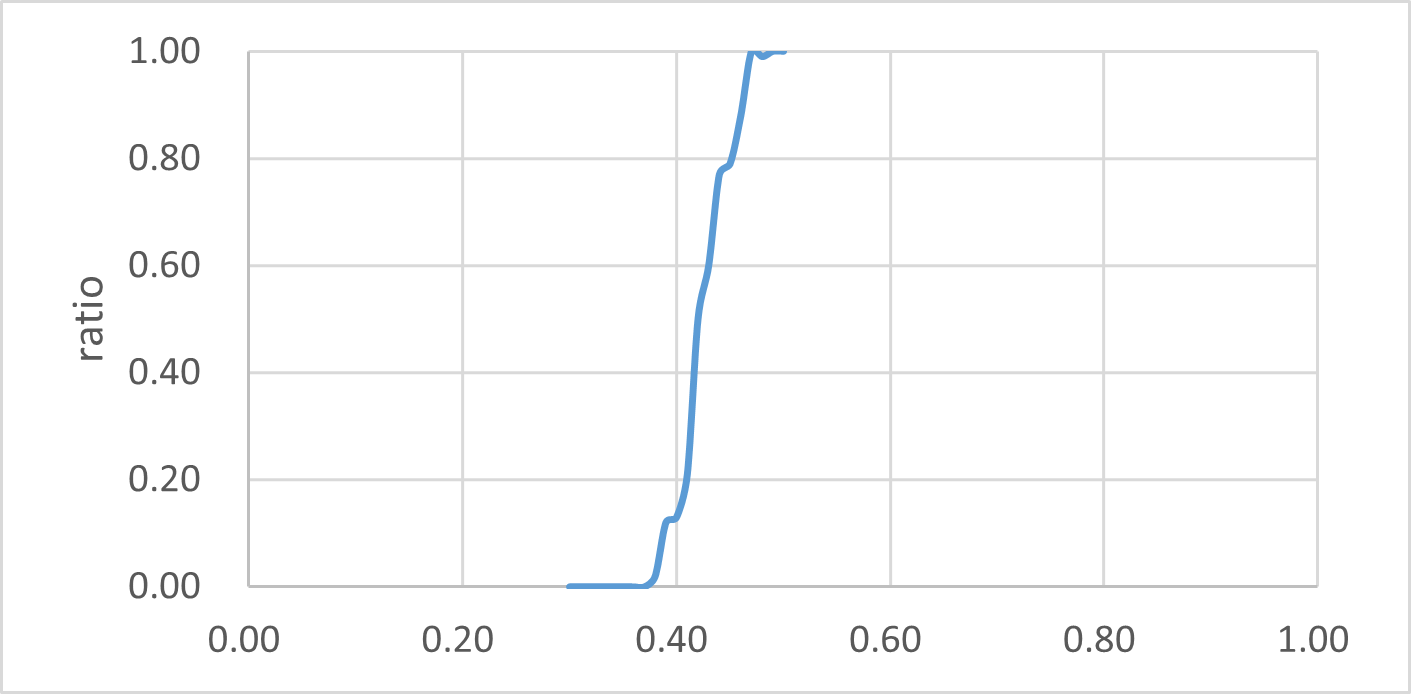}  
			\caption{m=200}
			\label{fig:(b)}
		\end{subfigure}
		\begin{subfigure}{0.45\textwidth}
			\centering
			\includegraphics[width=\textwidth]{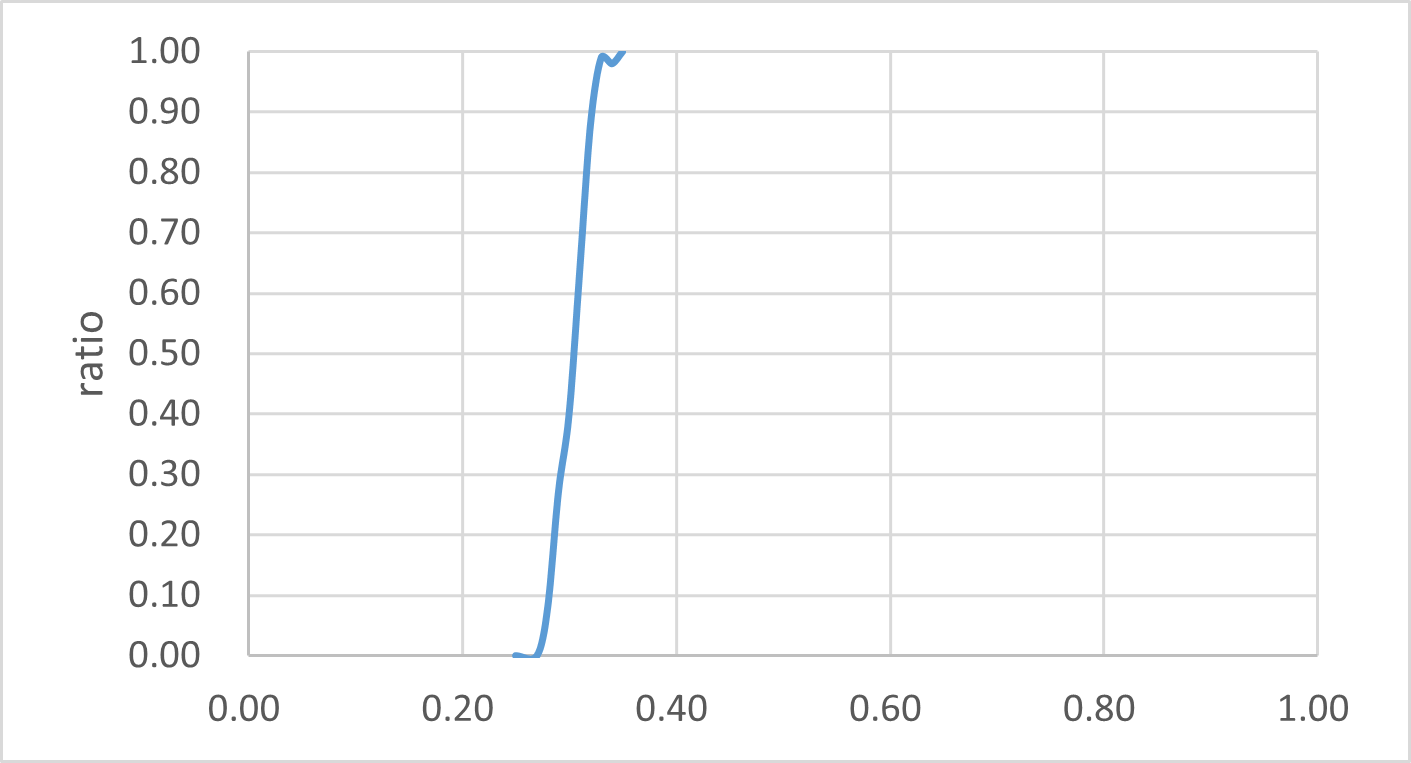}  
			\caption{m=400}
			\label{fig:(c)}
		\end{subfigure}
		\hfill
		\begin{subfigure}{0.45\textwidth}
			\centering
			\includegraphics[width=\textwidth]{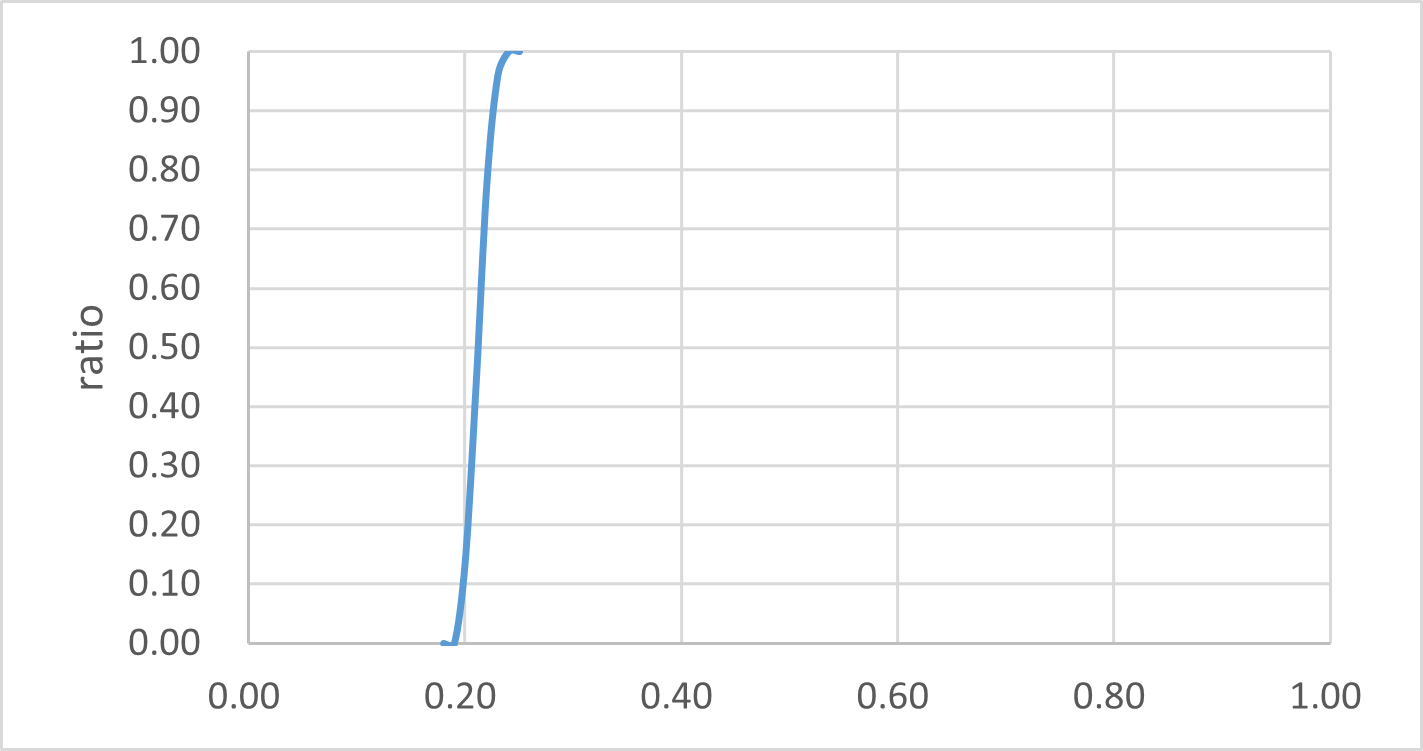}  
			\caption{m=800}
			\label{fig:(d)}
		\end{subfigure}
		
		\caption{The relationship between data density and recovery ratio across varying dimensions}
		\label{Fig1}
	\end{figure}

\begin{figure}[ht]
    \centering
    \begin{minipage}[b]{0.48\textwidth}
        \centering
        \includegraphics[width=\textwidth]{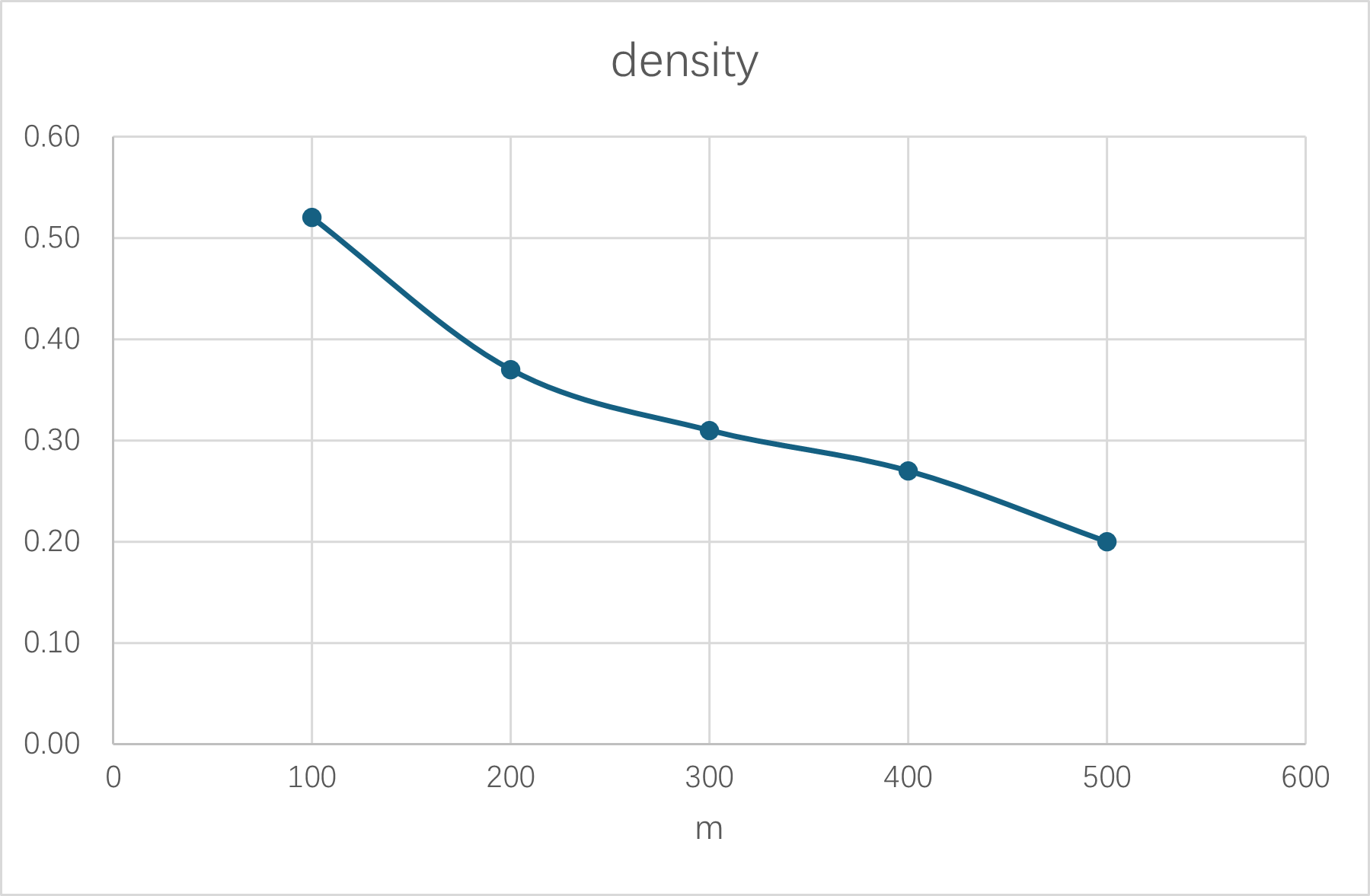}
        \caption{The effect of dimension on data density}
        \label{Fig2}
    \end{minipage}
    \hfill
    \begin{minipage}[b]{0.48\textwidth}
        \centering
        \includegraphics[width=\textwidth]{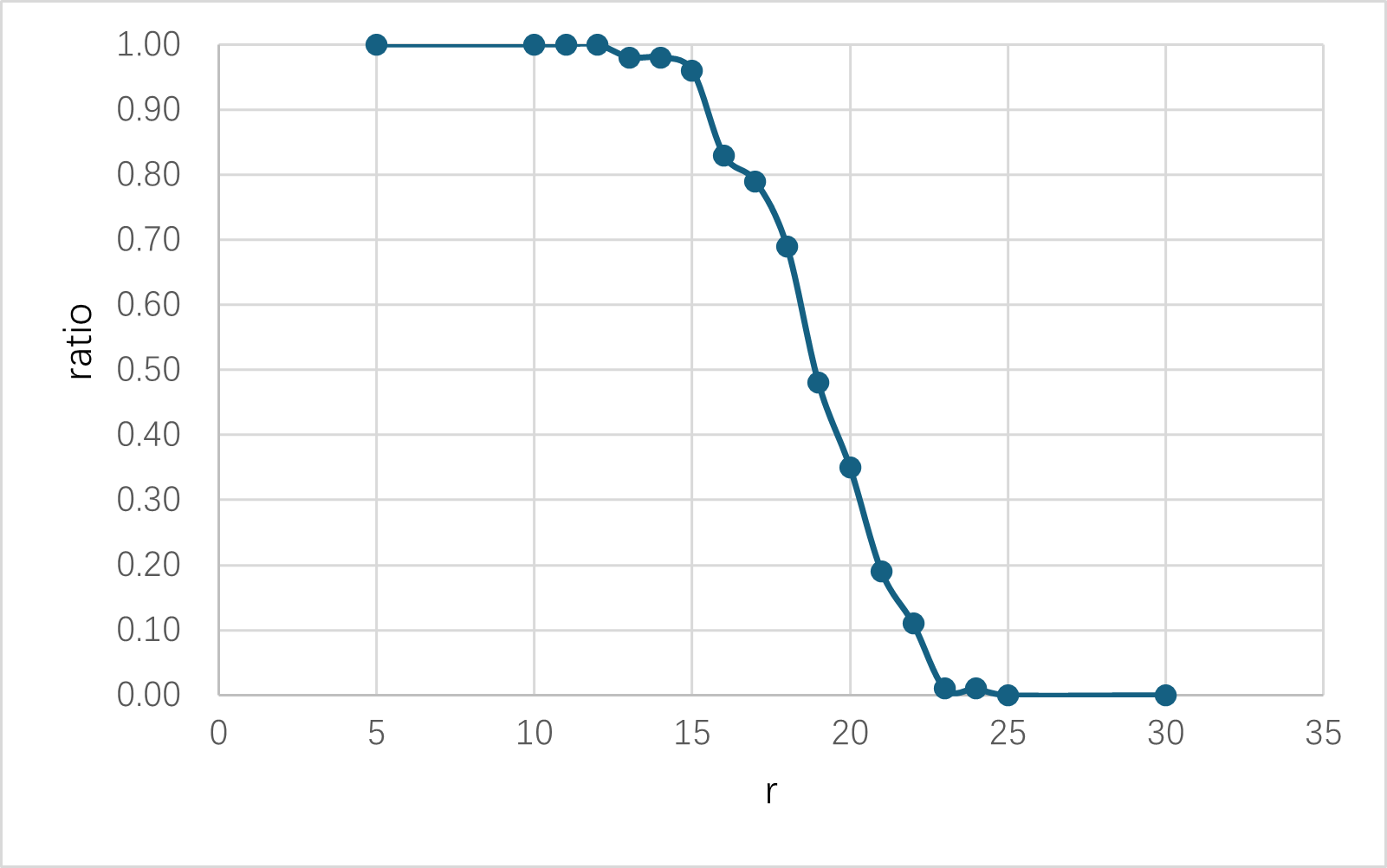}
        \caption{The relationship between rank $r$ and recovery ratio}
        \label{Fig3}
    \end{minipage}
\end{figure}

	We first investigate the impact of network scale $m$ and observation density on the ratio of exact recovery, which is defined as the proportion of independent trials in which the algorithm successfully reconstructs the true underlying low-rank matrix. For dimensions $m = 100, 200, 400, 800$, the exact recovery ratio increases with density (see Fig.~\ref{Fig1}); at the same density the critical density (the minimal observation density threshold required for the recovery ratio drops to 0) decreases with larger $m$ (see Fig.~\ref{Fig2}). The results indicate that lower density or larger matrix size makes exact recovery more difficult, consistent with the theoretical sensitivity of low-rank symmetric factorization under sparse masking.

    \subsubsection{Impact of Rank $r$ and Noise Level}
    
	We further fixed $m = 200$ and density $ = 0.7$ to observe the effect of the latent feature dimension $r$. The results indicate that an increase in $r$ weakens the uniqueness of the matrix factorization, thereby lowering the probability of exact recovery. 
	This phenomenon reflects degradation of the low-rank assumption at high rank and provides guidance for practical selection of community number. Additionally, we compared the impact of varying noise levels. Regardless of whether the network is masked or not, lower noise levels consistently yield higher recovery precision. This corroborates the theoretical reliability of our proposed model under ideal conditions, as well as its tolerance to moderate noise.
    		\begin{figure}[htbp]
		\centering
		\begin{subfigure}{0.45\textwidth}
			\centering
			\includegraphics[width=\textwidth]{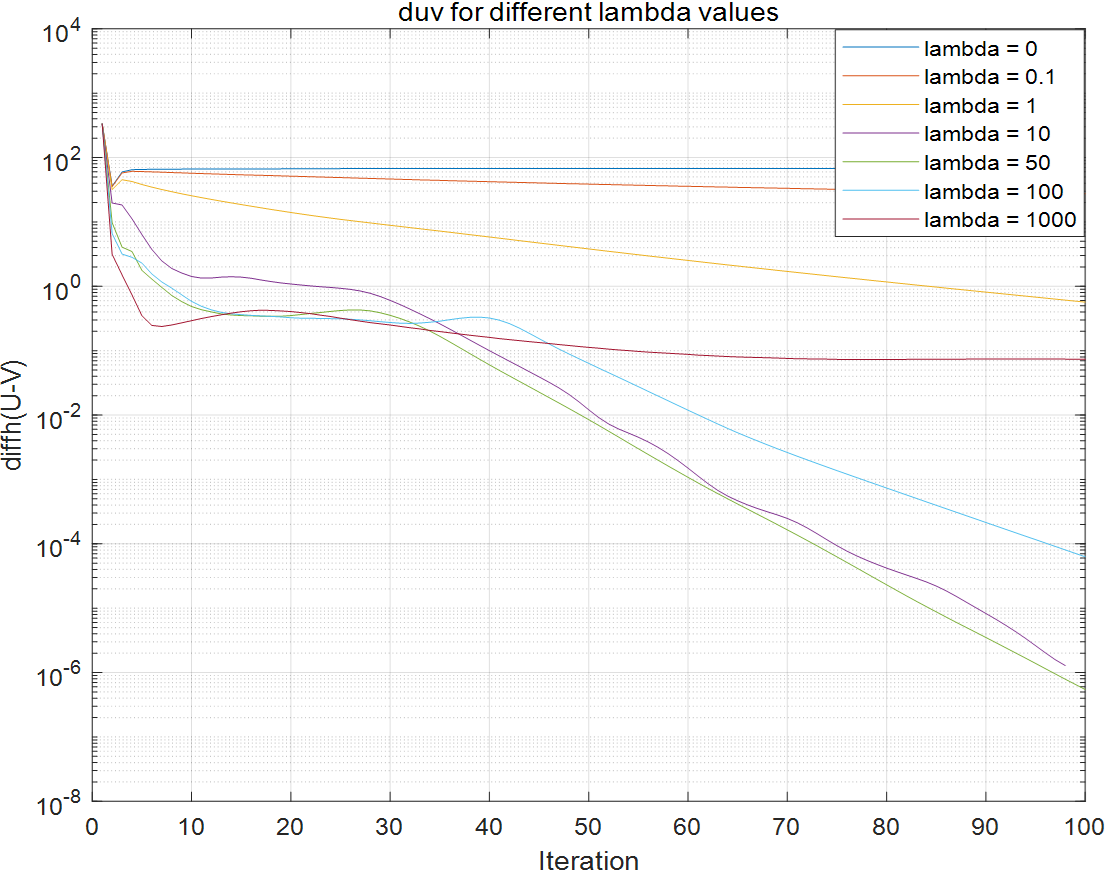}
			\caption{$\lambda$ to Diff$(U,V)$}
			\label{Fig4(a)}
		\end{subfigure}
		\begin{subfigure}{0.45\textwidth}
			\centering
			\includegraphics[width=\textwidth]{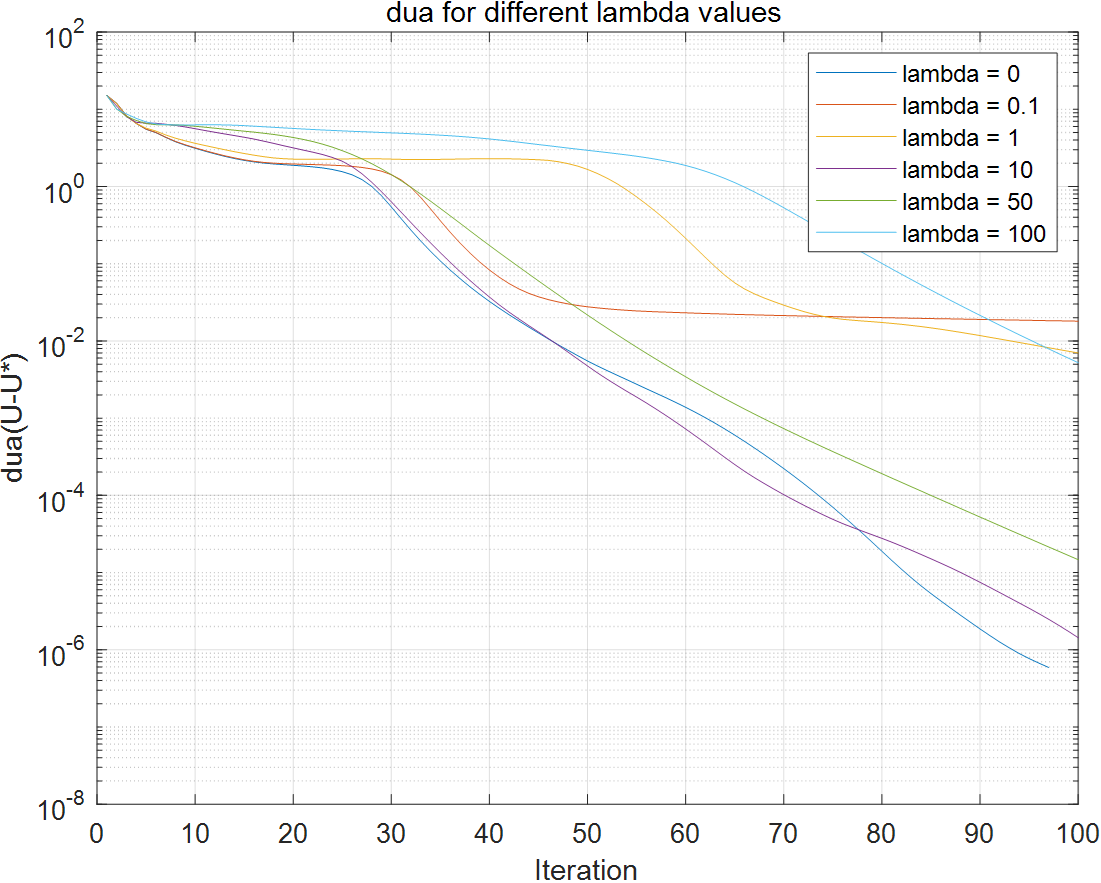}
			\caption{$\lambda$ to Diff$(U,U^*)$}
			\label{Fig4(b)}
		\end{subfigure}
		\caption{Comparison of algorithm performance under different $\lambda$ values}
		\label{Fig4}
	\end{figure}
    
	\subsubsection{Impact of Penalty Parameter $\lambda$}
	
	The regularization parameter $\lambda$ plays a critical role in balancing the approximation error and the symmetry constraint. We test the impact of various $\lambda$ values on the algorithmic convergence under combinations of noisy and noise-free, as well as masked and mask-free scenarios.
	
	We first examine the effect of $\lambda$ on feasible solutions. Figure~\ref{Fig4(a)} shows the variation in the difference between factor matrices $\text{Diff}(U, V)$ under different values of $\lambda$. The result indicates that the two factor matrices $U$ and $V$ eventually become identical given a sufficient number of iterations. In the presence of masks, this difference curve does not decrease monotonically as $\lambda$ increases, but exhibits a ``down-then-up before stabilizing'' oscillatory behavior.
	Furthermore, we investigate the effect of $\lambda$ on the rate of convergence. Figure~\ref{Fig4(b)} displays the distance to the true solution, $\text{Diff}(U, U^*)$. The results demonstrate that a larger $\lambda$ does not necessarily guarantee faster convergence; instead, regardless of the presence of noise or masks, the error curves undergo a turning phase where they first rise before falling. Synthesizing the tests, the algorithm exhibits the optimal convergence speed and solution accuracy when $\lambda$ is set between 10 and 50.
	These observations confirm the effectiveness of the exact penalty function under masking and suggest selecting $\lambda$ according to noise and missing conditions.

	
	Collectively, the above synthetic experiments clearly quantify the influence of each key parameter and furnish a reliable basis for parameter tuning and algorithm selection in real-network experiments.
	
	\subsection{Real-world Dataset}
	\subsubsection{Dataset}
	
	To further evaluate the Masked SymNMF model in a practical setting, we conduct experiments on a benchmark real‑world network: the Email‑Eu‑core dataset which is constructed from email communication data within a large European research institution and is publicly available from the Stanford Large Network Dataset Collection (SNAP).
	The network is defined as a directed graph, where a directed edge $u \to v$ exists if user $u$ has sent at least one email to user $v$. The dataset includes only internal communications among core members of the institution comprising a total of 1,005 nodes and 25,571 edges. Each node corresponds to a member and is associated with a ground‑truth community label indicating the department to which they belong. There are 42 departments in total, making the community detection task a 42‑class non‑overlapping partition problem.
	
	\subsubsection{Evaluation Protocol}
	
	Following~\cite{fortunato2010community,liu2023symmetry,cai2010graph,luo2020highly,vinh2009information,danon2005comparing}, we adopt Normalized Mutual Information (NMI) and Adjusted Rand Index (ARI) to evaluate the accuracy of involved community detectors. NMI measures the similarity between the resulting community division and the ground‑truth community information, while ARI assesses the agreement between two partitions corrected for chance. Both metrics are widely used in community detection literature and provide complementary perspectives on clustering quality.
	
	Let \(C = \{C_1, C_2, \dots, C_K\}\) denote the set of ground‑truth communities and \(R = \{R_1, R_2, \dots, R_K\}\) the community detection results obtained by a detector. The NMI between \(C\) and \(R\) is defined as:
	\begin{equation}
		\label{eq:50}
		\text{NMI}(C,R) = \frac{\displaystyle \sum_{g=1}^{K} \sum_{r=1}^{K} n_{g,r} \log \frac{n_{g,r} \cdot n}{n_{C,g} \cdot n_{R,r}}}{\sqrt{\left( \sum_{g=1}^{K} n_{C,g} \log \frac{n_{C,g}}{n} \right) \left( \sum_{r=1}^{K} n_{R,r} \log \frac{n_{R,r}}{n} \right)}},
	\end{equation}	
	where \(n\) is the total number of nodes, \(n_{g,r}\) is the number of nodes that belong to ground‑truth community \(C_g\) and are assigned to detected community \(R_r\), \(n_{C,g}\) is the number of nodes in \(C_g\), and \(n_{R,r}\) is the number of nodes in \(R_r\). NMI takes values in \([0,1]\), with larger values indicating higher agreement between the two partitions~\cite{strehl2002cluster}.
	
	The Adjusted Rand Index (ARI) measures the similarity between two partitions while correcting for the expected similarity due to random chance~\cite{hubert1985comparing,vinh2009information}. Given the contingency table where entry \(n_{g,r}\) denotes the number of nodes common to ground‑truth community \(C_g\) and detected community \(R_r\), the ARI is computed as:
	\begin{equation}
		\label{eq:51}
		\text{ARI}(C,R) = \frac{ \sum_{g,r} \binom{n_{g,r}}{2} - \left[ \sum_g \binom{n_{C,g}}{2} \sum_r \binom{n_{R,r}}{2} \right] \big/ \binom{n}{2} }{ \frac{1}{2} \left[ \sum_g \binom{n_{C,g}}{2} + \sum_r \binom{n_{R,r}}{2} \right] - \left[ \sum_g \binom{n_{C,g}}{2} \sum_r \binom{n_{R,r}}{2} \right] \big/ \binom{n}{2} }.
	\end{equation}
	The ARI ranges from \(-1\) to \(1\), with 1 indicating perfect agreement, 0 indicating random agreement, and negative values indicating agreement worse than random~\cite{hubert1985comparing,vinh2009information}.
	
	\subsubsection{Tested Models and Settings}
	For a comprehensive comparison, we evaluate the basic Masked SymNMF and its graph-regularized variants (\Cref{Section4}). Initial factors \(U_0,V_0\) are generated via singular value decomposition (SVD). Missing rates are set to \(\{0,0.2,0.4,0.6,0.8\}\), with all parameters properly tuned.
	
	To ensure convergence and efficiency, we adopt a dual stopping criterion: the projected gradient norm or objective change over 100 iterations drops below \(10^{-6}\). For PGD, an early stop over a 50-iteration window is used to mitigate oscillations, ensuring stable results.
	
	\subsubsection{Performance Verification}
    We evaluate our masked framework against traditional SymNMF on the Email-Eu-core dataset, using HALS, MU, and PGD. Traditional SymNMF imputes missing entries with zeros, while our Masked SymNMF uses a mask to fit only observed data. Both methods are tested across mask rates \(\{0.0,0.2,0.4,0.6,0.8\}\). We also validate graph-regularized Masked GSymNMF for further performance gains.
	
	\textit{A. Optimization and Convergence Analysis}
	From an optimization perspective, Masked SymNMF achieves lower final objective values than zero-imputed SymNMF (Table~\ref{tab:algorithm_comparison}), as fitting only observed entries via $\mathcal{P}_E$ avoids erroneous fitting to missing data and yields a tighter low-rank approximation.
    
	However, this comes with higher computational cost due to the masking operator and increased non-convexity under sparsity. Among the algorithms, PGD suffers severe oscillations and long runtimes at high mask rates, while HALS and MU converge more robustly. Masked HALS offers the best trade-off, balancing convergence speed and accuracy with manageable runtime.
	
	\begin{table}[htbp]
		\centering
        
		\caption{Performance comparison of different algorithms}
		\label{tab:algorithm_comparison}
		\begin{tabular}{@{}lccccc@{}}
			\toprule
			{Algorithm} & {Final Objective} & {$U-V$} & {Runtime} & {NMI} & {ARI} \\
			\midrule
			Masked HALS & 4.8736e+03 & 2.4425e+01 & 357.6606  & 0.5275 & 0.2683 \\
			HALS     & 9.8969e+03 & 1.4564e+01 & 1.2497    & 0.4640 & 0.1687 \\
			Masked MU   & 3.5196e+03 & 8.7721e+00 & 694.2600  & 0.5040 & 0.2330 \\
			MU       & 9.5516e+03 & 1.8483e+01 & 157.6402  & 0.4678 & 0.1942 \\
			Masked PGD   & 3.4094e+03 & 7.4934e+00 & 71315.9890 & 0.4814 & 0.1918 \\
			PGD       & 9.5426e+03 & 1.9616e+01 & 10.8154   & 0.4574 & 0.1931 \\
			\bottomrule
		\end{tabular}
	\end{table}
    
	\textit{B. Community Detection Performance}
	
\begin{figure}[htbp]
		\centering
		\includegraphics[width=0.65\textwidth]{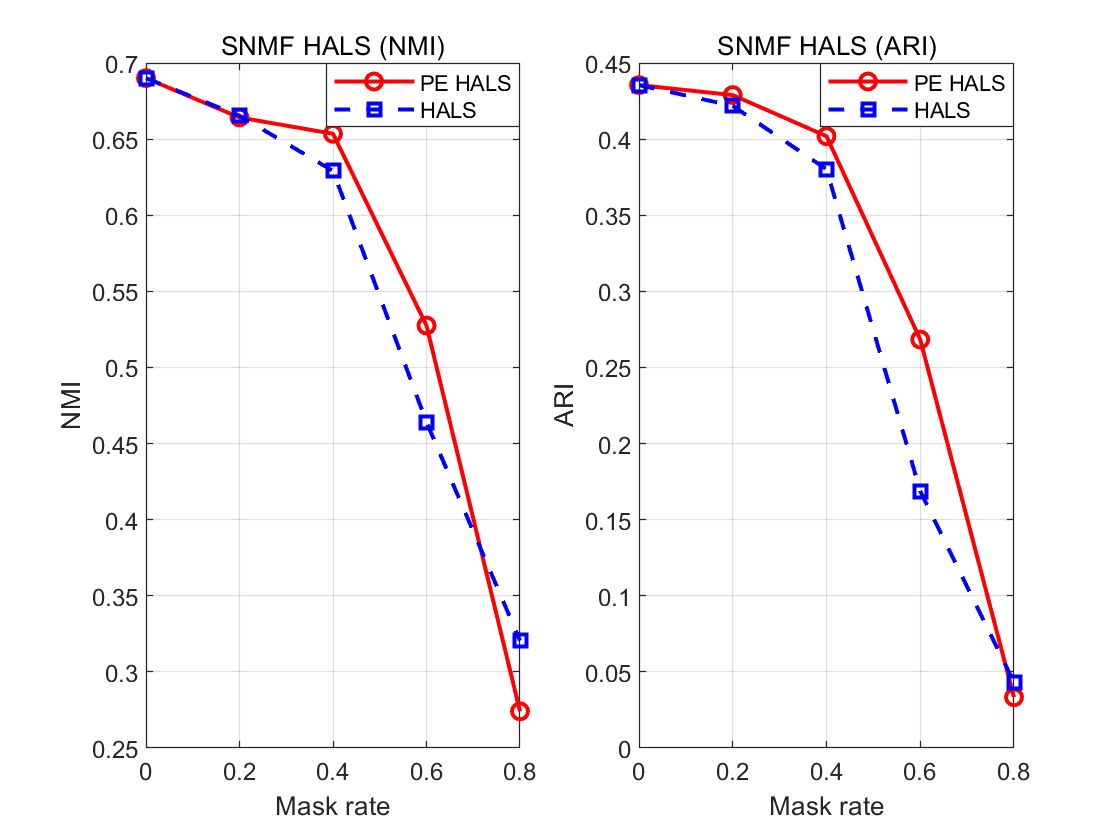}  
		\caption{NMI and ARI performance of the SymNMF HALS algorithm under different mask rates.}
		\label{Fig5}
	\end{figure}

	The ultimate quality of community detection is quantified using NMI and ARI. As illustrated in the performance curves (Fig.~\ref{Fig5}), at low mask rates, Masked SymNMF performs comparably to the baseline SymNMF. However, as the degree of missing information increases, the flaws of the zero-imputation method become starkly apparent. When the mask rate reaches 0.6, the NMI and ARI of the traditional HALS algorithm plummet to 0.4640 and 0.1687 ($\lambda = 1$), respectively, whereas Masked HALS maintains scores of 0.5275 and 0.2683, exhibiting a significant performance advantage. This comparison demonstrates that the traditional method distorts the true network topology. In contrast, mask-based factorization remains faithful to the valid data, extracting more accurate latent community partitions in incomplete networks.

	\textit{C. Enhancements via Graph Regularization}
	
	The results indicate that after incorporating the graph regularization term, the masked model still outperforms the zero-imputation model. Furthermore, the graph regularization yields substantial benefits for community detection under masks. Comparing the basic Masked HALS with the graph-regularized Masked HALS, under different mask rates, both the NMI and ARI of the model improved, demonstrating better performance in community detection.
	Graph regularization enforces smoothness over the feature manifold, ensuring that nodes with true underlying connections maintain similar representations in the low-dimensional space. This mechanism ``bridges'' the information gaps caused by the masking operator.
	\begin{figure}[htbp]
		\centering
		\includegraphics[width=0.65\textwidth]{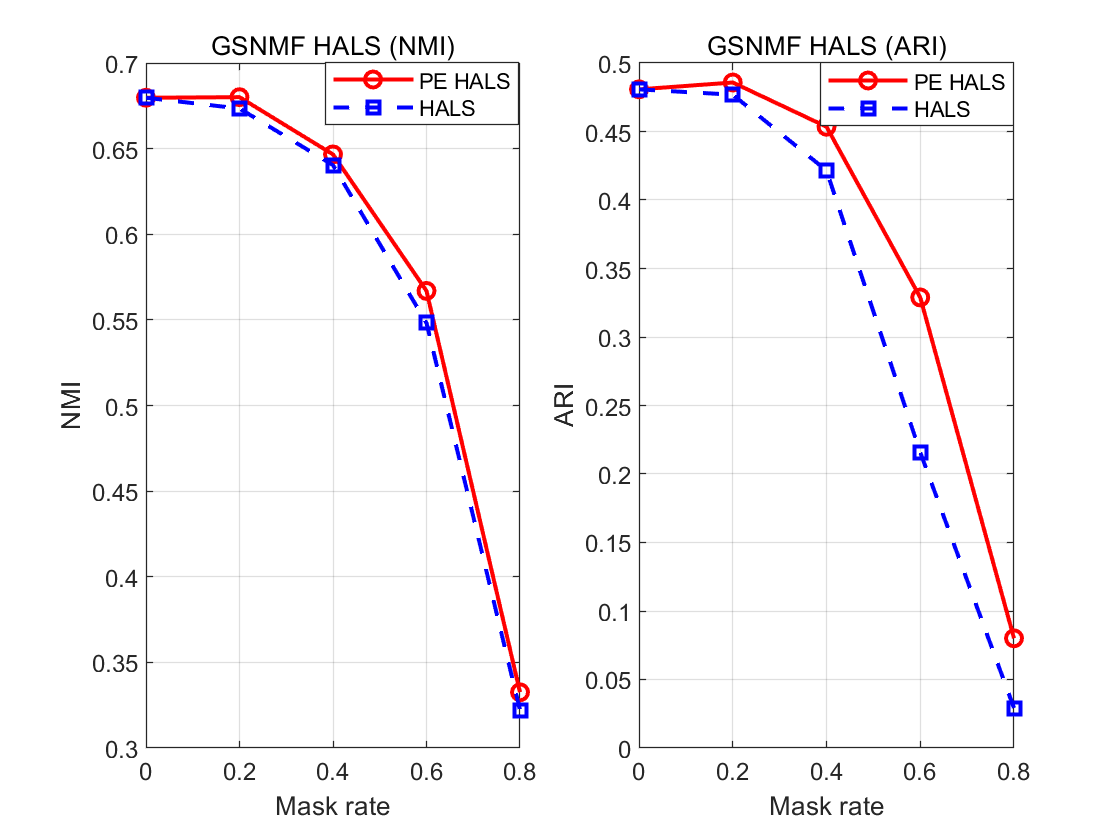}
		\caption{NMI and ARI performance of the GSymNMF HALS algorithm under different mask rates.}
		\label{Fig6}
	\end{figure}

    The experiments not only confirm the effectiveness of the proposed Masked SymNMF framework in handling sparse data but also demonstrate its extensibility. By integrating graph regularization (GSymNMF), the model can more fully exploit the limited observed structures. This confers practical value and competitiveness upon our approach when tackling the analysis of incomplete networks in the real world.
	\begin{table}[htbp]
		\centering
		\caption{Performance comparison of masked graph regularization model under different mask rates}
		\begin{tabular}{@{}lcccc@{}}
			\toprule
			Algorithm & Runtime & NMI & ARI & Mask rate \\
			\midrule
			Masked SymNMF HALS & 321.0561 & 0.6642 & 0.4289 & 0.2 \\
			Masked GSymNMF HALS & 1729.2936 & 0.6800 & 0.4858 & 0.2 \\
			Masked SymNMF HALS & 429.1247 & 0.6466 & 0.4019 & 0.4 \\
			Masked GSymNMF HALS & 2135.6827 & 0.6535 & 0.4535 & 0.4 \\
			Masked SymNMF HALS & 357.6606 & 0.5275 & 0.2683 & 0.6 \\
			Masked GSymNMF HALS & 1947.5323 & 0.5669 & 0.3289 & 0.6 \\
			Masked SymNMF HALS & 390.0678 & 0.2741 & 0.0333 & 0.8 \\
			Masked GSymNMF HALS & 797.8325 & 0.3325 & 0.0801 & 0.8 \\
			\bottomrule
		\end{tabular}
		\label{tab:alg_perf}
	\end{table}
	
	\section{Conclusion and Future Work}\label{Section6}
    We proposed Masked SymNMF for community detection in incomplete networks, using a masking operator and asymmetric relaxation with exact penalty. ANLS-based algorithms (MU, HALS, PGD) were developed. Experiments on synthetic and real data show consistent improvements over zero-imputation. Future work includes directed/dynamic networks, adaptive $\lambda$, and distributed implementations.

	\section*{Acknowledgement}
	R.G. was supported by the National Key R\&D Program of China under Grant No.~2022YFA1003800, the National Natural Science Foundation of China under Grant No.~12201318, and the Natural Science Foundation of Tianjin under Grant No.~25JCJQJC00300.
	R.J.Z. was supported by the National Natural Science Foundation of China (No. 1250012017), the Natural Science Foundation of Tianjin, China (No. 24JCQNJC01970), and the Fundamental Research Funds for the Central Universities (No. 050-63263071).

    \subsection*{Data and Code Availability}
    The source code for the proposed Masked SymNMF model and all experiments in this paper is available at \url{https://github.com/liuanqi7/Masked-SymNMF}.
    
	\bibliographystyle{siamplain}
	\bibliography{references}
\end{document}

%% file: ex_shared.tex
\usepackage{lipsum}
\usepackage{amsfonts}
\usepackage{graphicx}
\usepackage{epstopdf}
\usepackage{algpseudocode}
\usepackage{booktabs}
\usepackage{subcaption}
\usepackage{siunitx}
\usepackage{algorithm}
\usepackage{hyperref}
\usepackage[capitalize, nameinlink]{cleveref}
\ifpdf
  \DeclareGraphicsExtensions{.eps,.pdf,.png,.jpg}
\else
  \DeclareGraphicsExtensions{.eps}
\fi

\newsiamremark{remark}{Remark}
\newsiamremark{hypothesis}{Hypothesis}
\crefname{hypothesis}{Hypothesis}{Hypotheses}
\newsiamthm{claim}{Claim}
\newsiamremark{fact}{Fact}
\crefname{fact}{Fact}{Facts}

\headers{Masked SymNMF for Community Detection in Incomplete Networks}{A. Liu, R. Gu, and R. Zhang}

\title{Masked Symmetric Nonnegative Matrix Factorization for Community Detection in Incomplete Networks}

\author{Anqi Liu\thanks{School of Statistics and Data Science, Nankai University, Tianjin 300071, China
  (\email{2120240151@mail.nankai.edu.cn}). The first two authors contributed equally to this work.}
\and Ran Gu\footnotemark[1] \thanks{NITFID, School of Statistics and Data Science, LPMC, KLMDASR, and AAIS, Nankai University, Tianjin 300071, China
  (\email{rgu@nankai.edu.cn}).}
\and Rui-Jin Zhang\thanks{Corresponding author. School of Mathematical Sciences, Nankai University, Tianjin 300071, China
  (\email{zhangrj@nankai.edu.cn}).}}

\usepackage{amsopn}
